\documentclass[11pt]{article}
\usepackage[dvips]{epsfig}

\newtheorem{theorem}{Theorem}[section]
\newtheorem{proposition}[theorem]{Proposition}
\newtheorem{lemma}[theorem]{Lemma}
\newtheorem{corollary}[theorem]{Corollary}

\newcommand{\T}{{\cal T}}

\newcommand{\Qed}{\rule{2.5mm}{3mm}}

\newcommand{\mod}{\, \hbox{mod} \, }

\newcommand{\Aut}{\hbox{Aut}\, }

\newenvironment{proof}{{\noindent \sc Proof.}}{\hfill $\Qed$ \\}

\def\ZZ{{\hbox{\sf Z\kern-.43emZ}}}

\newcommand{\GP}{\hbox{GP}}

\begin{document}


\begin{center}
{\bf\large DISTANCE-BALANCED GRAPHS: SYMMETRY CONDITIONS}
\end{center}

\bigskip
\bigskip

\begin{center}
 Klavdija Kutnar,\footnotemark
 ~Aleksander Malni\v c,\footnotemark
 ~Dragan Maru\v si\v c\addtocounter{footnote}{-2} \footnotemark $^,$\footnotemark $^,$*
 ~and \v Stefko Miklavi\v c \addtocounter{footnote}{-2}\footnotemark
 $^,$\footnotemark

 \end{center}

\bigskip
\begin{center}
University of Primorska, Cankarjeva 5, 6000 Koper, Slovenia, and \\
Institute of Mathematics, Physics, and Mechanics, Jadranska 19, 1111 Ljubljana, Slovenia
\end{center}

\addtocounter{footnote}{-1} \footnotetext{~University of Primorska.
Supported in part by ``Ministrstvo za visoko \v solstvo, znanost in
tehnologijo Slovenije'', program no. P1-0285.}
\addtocounter{footnote}{1} \footnotetext{~Institute of Mathematics,
Physics, and Mechanics. Supported in part by ``Ministrstvo za visoko
\v solstvo, znanost in tehnologijo Slovenije'', program no. P1-0285.

 ~* Corresponding author e-mail: ~dragan.marusic@guest.arnes.si}

\vskip1.3cm
\begin{abstract}
A graph $X$ is said to be {\it distance--balanced} if for any edge
$uv$ of $X$, the number of vertices closer to $u$ than to $v$ is
equal to the number of vertices closer to $v$ than to $u$. A graph
$X$ is said to be {\it strongly distance--balanced} if for any edge
$uv$ of $X$ and any integer $k$, the number of vertices at distance
$k$ from $u$ and at distance $k+1$ from $v$ is equal to the number
of vertices at distance $k+1$ from $u$ and at distance $k$ from $v$.
Obviously, being distance--balanced is metrically a weaker condition
than being strongly distance--balanced. In this paper, a connection
between symmetry properties of graphs and the metric property of
being (strongly) distance--balanced is explored. In particular, it
is proved that every vertex--transitive graph is strongly
distance--balanced.

A graph is said to be {\em semisymmetric} if its automorphism group
acts transitively on its edge set, but does not act transitively on
its vertex set. An infinite family of semisymmetric graphs, which
are not distance--balanced, is constructed.

Finally, we give a complete classification of strongly
distance--balanced graphs for the following infinite families of
generalized Petersen graphs: $\GP(n,2)$, $\GP(5k+1,k)$, $\GP(3k\pm
3,k)$, and $\GP(2k+2,k)$.

\end{abstract}


\section{Introduction}
\label{sec:intro}
\indent

Let $X$ be a graph with diameter $d$, and let  $V(X)$ and $E(X)$
denote the vertex set and the edge set of $X$, respectively. For
$u,v \in V(X)$, we let $d(u,v)$ denote the minimal path-length
distance between $u$ and $v$. We say that $X$ is {\em
distance--balanced} if
$$
  |\{x\in V(X)\mid d(x,u)<d(x,v)\}| = |\{x\in V(X)\mid d(x,v)<d(x,u)\}|
$$
holds for an arbitrary pair of adjacent vertices $u$ and $v$ of $X$.
These graphs were,  at least implicitly,  first studied by
Handa \cite{Ha} who considered distance--balanced partial cubes.
The term itself,  however,  is due to Jerebic, Klav\v zar and Rall \cite{JKR}
who studied  distance--balanced graphs in the framework of various kinds of graph products.

\begin{figure}[ht]
\begin{center}
 \includegraphics[width=0.70\hsize]{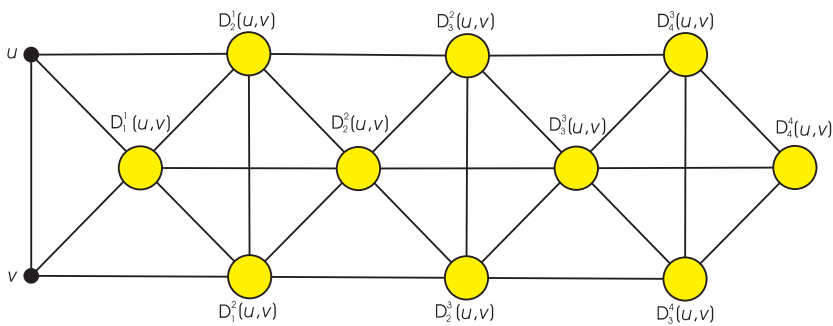}
\caption{\small A distance partition of a graph with diameter $4$
with respect to edge $uv$. 
\label{struc}}
\end{center}
\end{figure}

Let  $uv$ be an arbitrary edge of $X$. For any two integers $k,l$,
we
let
$$
  D^k_l(u,v)=\{x \in V(X)\mid d(u,x)=k \ \ \ \textrm{and} \ \ \ d(v,x)=l\}.
$$
The triangle inequality implies that only the sets $D^{k-1}_k(u,v)$,
$D^{k}_k(u,v)$ and $D^{k}_{k-1}(u,v)$ for $k \in \{1, \ldots, d\}$
can be nonempty. The sets $D^k_l(u,v)$ give rise to  a "distance
partition" of $V(X)$ with respect to the edge $uv$ (see
Figure~\ref{struc}). Moreover, one can easily see that $X$ is
distance--balanced if and only if

\begin{equation}
\label{equ:1}
  \sum_{k=1}^{d}|D^{k}_{k-1}(u,v)|=\sum_{k=1}^{d}|D^{k-1}_k(u,v)|
\end{equation}
holds for every edge $uv\in E(X)$.

 \medskip

\medskip

Obviously,  if $|D^k_{k-1}(u,v)|=|D^{k-1}_k(u,v)|$ holds for
$1\le k\le d$ and for every edge $uv\in E(X)$, then $X$ is
distance--balanced. The converse, however,  is not necessarily true.
For instance, in the generalized Petersen
graphs $\GP(24,4)$, $\GP(35,8)$ and $\GP(35,13)$
(see Section \ref{sec:strong} for a formal definition), we can find
two adjacent vertices $u$, $v$ and an integer $k$,
such that $|D^k_{k-1}(u,v)| \ne |D^{k-1}_k(u,v)|$. But it is
easy to see that these graphs are distance--balanced.

We shall say that $X$ is  {\it strongly distance--balanced}, if
$|D^{k}_{k-1}(u,v)|=|D^{k-1}_k(u,v)|$ for every integer $k$
and every edge $uv\in E(X)$.
Observe that distance--regular graphs are strongly distance--balanced.
(We refer the reader to  \cite{BCN} for the
definition and basic properties of distance--regular graphs.)
Being strongly distance--balanced is therefore metrically a weaker
condition than being distance--regular.
It is well known that not every
distance--regular graph is vertex--transitive (see \cite[p. 139]{NB74} for an example),
and thus not every distance--balanced graph is vertex--transitive.

The aim of this article is to explore a purely metric property of
being (strongly) distance--balanced in the context of graphs
enjoying certain special symmetry conditions. For example,  as
observed in \cite[Prop. 2.4]{JKR}, it is obvious that
arc--transitive graphs are necessarily distance--balanced. Namely,
such graphs contain automorphisms which interchange adjacent
vertices. A general vertex--transitive graph, however, may contain
edges which are not flipped over by an automorphism and therefore it
is not immediately obvious that it should be distance--balanced. But
as we shall see in Corollary~\ref{cor:vertex}, vertex--transitive
graphs are not only distance--balanced; they are also strongly
distance--balanced. Furthermore,  since being vertex--transitive is
not a necessary condition for a graph to be distance--balanced,
studying graphs which are as close to being vertex--transitive as
possible, seems like the next step to be taken. In Section
\ref{sec:semisym} we construct an infinite family of
edge--transitive but not vertex--transitive graphs which are not
distance--balanced (see Proposition~\ref{prop:folk}).

Finally, in Section \ref{sec:strong}, we explore the property of
being strongly distance--balanced for the family of generalized
Petersen graphs. We give a complete classification of strongly
distance--balanced graphs for the following infinite families:
 $\GP(n,2)$ for $n \ge 3$ and $n \ne 4$ (see Proposition \ref{n2}),
$\GP(5k+1,k)$ (see Proposition \ref{5k+1}),
$\GP(3k\pm 3,k)$ (see Theorem \ref{thm:glavni}),
and $\GP(2k+2,k)$ (see Theorem \ref{thm:2k+2}).


\section{Vertex--transitive graphs}
\label{sec:basic} \indent

In this section we give a characterization of  strongly
distance--balanced graphs, and as a consequence prove that every
vertex-transitive graph is strongly distance--balanced. Recall that
a graph  $X$,  with vertex set $V(X)$, edge set $E(X)$, arc set
$A(X)$ and the automorphism group $\Aut X$, is said to be {\em
vertex--transitive}, {\em edge--transitive}, and {\em
arc--transitive}, if $\Aut X$ acts transitively on $V(X)$, $E(X)$,
and $A(X)$, respectively.

For a graph $X$, a vertex $u$ of $X$ and an integer $i$, let $S_i(u)
= \{x\in V(X)\mid d(x,u)=i\}$ denote the set of vertices of $X$
which are at distance $i$ from $u$. Let $u,v\in V(X)$ be adjacent
vertices. Observe that $S_i(u)$ is a disjoint union of the sets
$D^{i}_{i-1}(u,v),$ $D^{i}_{i}(u,v)$ and $D^{i}_{i+1}(u,v)$.
Similarly, $S_i(v)$ is a disjoint union of the sets
$D^{i-1}_i(u,v),$ $D^i_i(u,v)$ and $D^{i+1}_i(u,v)$.

\begin{proposition}
\label{pro:sphere} Let $X$  be a graph with diameter $d$. Then $X$
is strongly distance--balanced if and only if $|S_i(u)|=|S_i(v)|$
holds for every edge $uv \in E(X)$ and every $i \in \{0, \ldots, d\}$.
\end{proposition}

\begin{proof} Assume first that $X$ is strongly
distance--balanced and let $uv\in E(X)$. By definition, we have
$|D^i_{i+1}(u,v)|=|D^{i+1}_i(u,v)|$ for $i \in \{0, \ldots, d-1\}$. However, since
$S_i(u)=D^{i}_{i-1}(u,v)\cup D^{i}_{i}(u,v)\cup D^{i}_{i+1}(u,v)$ (disjoint union), and
$S_i(v)=D^{i-1}_{i}(u,v)\cup D^{i}_{i}(u,v)\cup D^{i+1}_{i}(u,v)$ (disjoint union),
we have also $|S_i(u)|=|S_i(v)|$ for $i \in \{0, \ldots, d\}$.

\medskip

Next assume that $|S_i(u)|=|S_i(v)|$ holds for every edge $uv$ of
$X$ and every $i \in \{0, \ldots, d\}$. Using induction we now show
that $|D^i_{i+1}(u,v)|=|D^{i+1}_i(u,v)|$ holds for every edge $uv$
of $X$ and every $i  \in \{0, \dots, d-1\}$. Obviously,
$|D^0_1(u,v)|=|D^1_0(u,v)|=1$. Suppose now that
$|D^{k-1}_k(u,v)|=|D^k_{k-1}(u,v)|$ for some $1 \le k\le d-1$.
Observe that
$$
  |D^k_{k+1}(u,v)|=|S_k(u)|-|D^k_k(u,v)|-|D^k_{k-1}(u,v)|
$$
and
$$
  |D^{k+1}_k(u,v)|=|S_k(v)|-|D^k_k(u,v)|-|D^{k-1}_k(u,v)|.
$$
Since $|S_k(u)|=|S_k(v)|$ and in view of the induction hypothesis also
$|D^{k-1}_k(u,v)|=|D^k_{k-1}(u,v)|$, we obtain
$$
  |D^k_{k+1}(u,v)|=|D^{k+1}_k(u,v)|.
$$
The result follows.
\end{proof}

Let $X$ be a connected strongly distance--balanced graph with diameter $d$.
Then, by Proposition \ref{pro:sphere}, $|S_i(u)|=|S_i(v)|$
holds for any pair of adjacent vertices $u,v\in V(X)$ and every
$i  \in \{0, \dots, d\}$.
Observe that connectedness implies that $|S_i(u)|=|S_i(v)|$
holds for any pair of vertices $u,v\in V(X)$ and every
$i  \in \{0, \dots, d\}$.
Let us remark that graphs with this property are also called
{\em distance--degree regular}. Distance--degree regular graphs were
studied in \cite{HN84}.

\smallskip
Since automorphisms preserve distances, we have the following
immediate consequence for vertex--transitive graphs.

\begin{corollary}
\label{cor:vertex} Every vertex--transitive  graph is
strongly distance--balanced.
\end{corollary}


\section{Semisymmetric graphs}
\label{sec:semisym}
\indent

\begin{figure}[ht!]
\begin{footnotesize}
\begin{center}
\includegraphics[width=0.30\hsize]{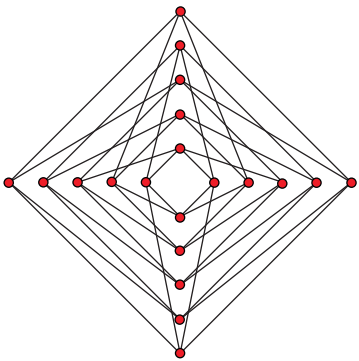}
\caption{\label{Folkman}\footnotesize The Folkman graph.}
\end{center}
\end{footnotesize}
\end{figure}

A regular edge--transitive graph which is not vertex--transitive is
usually called {\em semisymmetric}. Note that a semisymmetric graph
is necessarily bipartite, where the two sets of bipartition coincide
with the two orbits of the automorphism group. The smallest
semisymmetric graph has 20 vertices and was discovered by Folkman
\cite{F67} (see Figure~\ref{Folkman}) who initiated this topic of
research. Since then the theory of semisymmetric graphs has come a
long way (see
\cite{CMMP,DX99,I87,LX02,LXW04,MMPW04,MP01,MP02,Park05}).

\bigskip

\bigskip

As we have seen in the previous section,
vertex--transitive graphs are distance--balanced.
It is therefore natural to explore
the property of being distance--balanced within
the class of semisymmetric graphs; a class of objects which
is as close to vertex--transitive graphs as one can possibly get.

Given a graph $X$, we may associate with each arc $x=(u,v)$ of $X$
the triple $(x_l,x_c,x_r)$, where $x_l$ is the cardinality of the
set of all those vertices which are closer to $u$ than they are to
$v$. Similarly, $x_r$ is the cardinality of the set of all those
vertices which are closer to $v$ than they are to $u$. And
finally, $x_c$ is the cardinality of the set of all those vertices which are
at equal distance from both $u$ and $v$. We call this triple the
{\em distance--balance triple} of the arc $x$.  Clearly, $X$ is
distance--balanced if $x_l=x_r$ for all arcs $x$ of $X$. Also, in
the case when $X$ is edge--transitive the distance--balance
triples are unique up to switching of left and right components.

Semisymmetric graphs have no automorphisms
which switch adjacent vertices,
and therefore, may arguably  be considered as good candidates
for graphs which are not distance--balanced.
However, there are semisymmetric graphs which are distance-balanced.
For example the Gray graph, the smallest cubic semisymmetric graph,
denoted by SS54 in \cite{CMMP}, is indeed not distance-balanced.
Its distance--balance triple is $(23,0,31)$.
(Note that the central component in a bipartite graph is always $0$.)
On the other hand, this triple is $(55,0,55)$ for the next smallest
cubic semisymmetric graph SS110 on $110$ vertices, and so this graph is
distance--balanced.

The object of this section is to present an infinite
family of semisymmetric graphs which are not distance--balanced.
The smallest member of this family is the Folkman graph mentioned above.

\bigskip

\begin{figure}[ht]
\begin{center}
 \includegraphics[width=0.25\hsize]{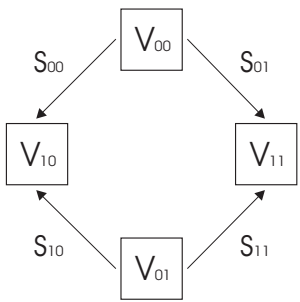}
\caption{\footnotesize Bipartite tetracirculant
$\T(n,S_{00},S_{01},S_{10},S_{11})$.\label{F-graph}}
\end{center}
\end{figure}

\bigskip

The construction comes from \cite{MP01}. We basically use the same
notation. Let $n$ be a positive integer and let $S_{00}$, $S_{01}$,
$S_{10}$ and $S_{11}$ be nonempty subsets of $\ZZ_n$. Define the
graph $X=\T(n,S_{00},S_{01},S_{10},S_{11})$ to have vertex set
$\ZZ_n\times \ZZ_2\times \ZZ_2$ and edge set $\{(a,0,i)(b,1,j)\mid
i,j\in \ZZ_2, b-a\in S_{ij} \}$. (The symbol $\T$ stand for
tetracirculant, a graph having an automorphism with four orbits of
equal length.) We use the shorthand notations $V_{00}=V_{00}(X)$,
$V_{11}=V_{11}(X)$, $V_{01}=V_{01}(X)$, and $V_{10}=V_{10}(X)$,
where $V_{ij}(X)=\{ (a,i,j) \mid a\in \ZZ_n\}$, $i,j\in \ZZ_2$.
Furthermore, we use the symbols $x_l$, $y_l$, $u_l$ and $w_l$, where
$l \in \{0,\dots,n-1\}$, to denote the elements of $V_{00}$,
$V_{11}$, $V_{01}$ and $V_{10}$, respectively (see
Figure~\ref{F-graph}). In particular, any graph of the form
$\T(n,R,R,T,T)$, where $R,T \subseteq \ZZ_n$, is called a {\em
generalized Folkman graph}. Let $p$ be a prime, let  $a\in \ZZ_p^*$,
and let $S$ be a nontrivial subgroup of $\ZZ_p^*$ such that
$a\not\in S$, but $a^2\in S$, and moreover $S\not=x+aS$ for all
$x\in \ZZ_p$, and $S\not=x+S$ for all $x\in \ZZ_p\setminus \{0\}$ .
Then $\T(p,S,S,aS,aS)$ is semisymmetric (see \cite{MP01}). (Here
$x+aS$ and $x+S$ are the sets $\{x+as\mid s\in S\}$ and $\{x+s\mid
s\in S\}$, respectively.)

In the special case when $S$ is the subgroup of all
squares in $\ZZ_p^*$, we use the symbol $N$ for the coset
$aS \neq S$ of all nonsquares. In this case the graphs $\T(p,S,S,aS,aS)$
have diameter equal to $4$. The smallest graph of this type is the above mentioned
Folkman graph $\T(5,S,S,N,N)$ with $20$ vertices, where $S=\{-1,1\}$
and $N=\{-2,2\}$ (see Figure~\ref{Folkman}).
These graphs are not distance-balanced, as is shown in the proposition below.

\begin{proposition}
\label{prop:folk} Let $p\geq5$ be a prime, and let $S$ and $N$ be
the set of squares and nonsquares in $\ZZ_p^*$. Then the
generalized Folkman graph $\T(p,S,S,N,N)$ is not
distance--balanced.
\end{proposition}

\begin{proof} Let $X=\T(p,S,S,N,N)$.
It is easy to see that  the Folkman graph ($p=5$)
is not distance--balanced.
We may therefore assume that $p>5$.
Since $X$ is regular and of diameter $4$
we have that, in view of (\ref{equ:1}),
it is sufficient to show that there exists an edge $uv \in E(X)$ such that

$$
  \sum_{k=1}^{3}D^k_{k+1}(u,v)\not=\sum_{k=1}^{3}D^{k+1}_k(u,v).
$$
Since $1\in S$, there exists an edge in
$X$ between $x_0\in V_{00}$ and $y_1\in V_{11}$.
It may be seen that
$$
  S_1(x_0)=\{y_s \mid s\in S\}\cup \{w_s\mid s\in S\},  \;
  S_2(x_0)=\{x_{i} \mid i\in \ZZ_p^*\}\cup  \{u_i \mid i\in \ZZ_p^*\},
$$
$$
  S_3(x_0) =  \{y_i\mid i\in \ZZ_p\setminus S\}\cup \{w_i \mid i\in \ZZ_p\setminus S\}, \;
  S_4(x_0) =  \{u_0\},
$$
and that
$$
  S_1(y_1)=\{x_{-s+1} \mid s\in S\}\cup \{u_{-as+1} \mid s\in S\},  \;
  S_2(y_1)=\{y_i \mid i\in \ZZ_p \setminus \{1\}\}\cup \{w_i \mid i \in \ZZ_p\},
$$
$$
  S_3(y_1)=\{x_{i} \mid i\in \ZZ_p\setminus \{1-s\mid s\in S\}\}\cup
  \{u_{i} \mid i\in \ZZ_p\setminus \{1-as\mid s\in S\}\}, \;
  S_4(y_1)=\emptyset.
$$
It follows that $|D^1_2(x_0,y_1)| =p-2$, $|D^2_3(x_0,y_1)| =p$
and $|D^{3}_{4}(x_0,y_1) |=0$. On the other hand, $|D^2_1(x_0,y_1)| =p-2$,
$|D^3_2(x_0,y_1)|=p+1$ and $|D^4_3(x_0,y_1) |=1$. Thus,
$$
  \sum_{k=1}^{3}|D^k_{k+1}(x_0,y_1)|\not=\sum_{k=1}^{3}|D^{k+1}_k(x_0,y_1)|.
$$
Therefore, $X$ is not distance--balanced.
\end{proof}


\section{Generalized Petersen graphs}
\label{sec:strong}
\indent

Let $n \ge 3$ be a positive integer, and let
$k \in \{1,\dots,n-1\} \setminus\{n/2\}$.
The generalized Petersen graph $\GP(n, k)$
is defined to have the following vertex set and edge set:
\begin{eqnarray}
\label{genpet}
  V(\GP(n,k)) &=& \{ u_i \mid i \in \ZZ_n \}\cup \{v_i \mid i \in \ZZ_n \},  \nonumber \\
  E(\GP(n, k)) &=& \{u_iu_{i+1} \mid i \in \ZZ_n\}\cup
    \{ v_{i}v_{i+k}  \mid i \in \ZZ_n\}\cup
    \{  u_iv_i \mid i \in \ZZ_n\}.
\end{eqnarray}

Note that $\GP(n, k)$ is cubic, and
that it is bipartite precisely when $n$ is even and $k$ is odd.
It is easy to see that $\GP(n,k)\cong\GP(n,n-k)$. Furthermore, if
the multiplicative inverse $k^{-1}$ of $k$ exists in $\ZZ_n$, then
the mapping $f$ defined by the rule

\begin{equation}
\label{iso}
f(u_i) = v_{k^{-1}i}, \qquad f(v_i) = u_{k^{-1}i}
\end{equation}

\noindent
gives us an isomorphism of graphs $\GP(n,k)$ and $\GP(n, k^{-1})$,
where the use of the same symbols for vertices in
$\GP(n,k)$ and $\GP(n, k^{-1})$ should cause no confusion.

In Section~\ref{sec:intro} we mentioned that not every
distance--balanced graph is also strongly distance--balanced.
Using program package {\sc Magma} \cite{Magma} one may easily see that
$\GP(24,4)$, $\GP(35,8)$ and $\GP(35,13)$  are
the only graphs among distance--balanced generalized Petersen graphs $\GP(n,k)$
on up to $120$ vertices which are not strongly distance--balanced.
This section is devoted to a more detailed investigation
of the property of being strongly
distance--balanced for several infinite families of the
generalized Petersen graphs.
We start with a rather straightforward observation.

\begin{proposition}
\label{n2}
Let $n \ge 3$ be an integer, $n \ne 4$. Then $\GP(n,2) $
is strongly distance--balanced if and only if $n \in \{3,5,7,10\}$.
\end{proposition}
\begin{proof}
It is easy to see that $|S_3(u_0)| = 6$ and $|S_3(v_0)|=4$ for $n \ge 13$.
Furthermore, if $n \le 12$ then $\GP(n,2)$ is strongly distance--balanced if and only if
$n \in \{3,5,7,10\}$.
\end{proof}

\medskip \noindent
The next proposition gives yet another infinite
family of generalized Petersen graphs for which it is easy to identify their
strongly distance--balanced members.

\begin{proposition}
\label{5k+1}
Let $k$ denote a positive integer. Then $\GP(5k+1,k)$ is strongly distance--balanced
if and only if $k=1$.
\end{proposition}
\begin{proof}
It can be easily verified that  $\GP(6,1)$ is the only
strongly distance--balanced graph for $k \le 5$.
As for $k \ge 6$, we have $|S_4(u_0)|=18$ and $|S_4(v_0)|=16$, and the result follows.
\end{proof}

\medskip
In order to investigate the property of being strongly distance--balanced
for certain other families of generalized Petersen graphs,
let us recall that the automorphism groups of the generalized Petersen
graphs were determinated in \cite{Frucht}.
Let $\rho, \tau: V(\GP(n,k)) \to V(\GP(n,k))$ be the mappings defined
by the rules $\rho(u_i) = u_{i+1}$, $\rho(v_i) = v_{i+1}$, $\tau(u_i) = u_{-i}$ and
$\tau(v_i) = v_{-i} \; (i \in \ZZ_n)$. Then

\begin{equation}
\label{grupa}
\langle \rho, \tau \rangle \subseteq \Aut(\GP(n, k)).
\end{equation}

\noindent
Moreover, $\GP(n,k)$ is vertex--transitive
if and only if $k^2 \equiv \pm 1 (\mod n)$ \cite{Frucht}.

Let us now analyze the family $\GP(3k+3,k)$, $k \ge 1$. To keep
things simple we assume that $k\ge 13$.

\begin{lemma}
\label{lem:k}
Let $k\ge 13$ be an integer, let $n=3k+3$, let $b=\lceil
(k+1)/2\rceil$, and let $u_0 \in V(\GP(n,k))$.
Then the following hold:
\begin{itemize}
\itemsep=0cm
\item[(i)] $S_1(u_0)=\{u_{\pm 1},v_0\}$, $S_2(u_0)=\{u_{\pm 2},v_{\pm 1},v_{\pm k}\}$,\\
          $S_3(u_0)=\{u_{\pm 3},u_{\pm k},v_{\pm 2},v_{\pm (k+1)}, v_{\pm (k-1)}, v_{\pm 2k}\}$,\\
          $S_4(u_0)=\{u_{\pm 4},u_{\pm (k+1)}, u_{\pm (k-1)},u_{\pm 2k}, v_{\pm 3},v_{\pm (k+2)}, v_{\pm (k-2)}, v_{\pm (k+4)}\}$,\\
          $S_5(u_0)=\{u_{\pm 5},u_{\pm (k+2)}, u_{\pm (k-2)},u_{\pm (k+4)}, v_{\pm 4},v_{\pm (k-3)},   v_{\pm (k+5)}\}$;

 \item[(ii)] $S_{i}(u_0)=\{u_{\pm i},
u_{\pm(k+i-1)}, u_{\pm(k-i+3)}, v_{\pm(i-1)}, v_{\pm( k-i+2)},v_{\pm(k+i)}\}$ for $6\le i\le b$;

 \item[(iii)] if $k$ is odd, then
$S_{b+1}(u_0)=\{u_{\pm((k+3)/2)}, v_{\pm((k+1)/2)},
u_{\pm((3k+1)/2)}, v_{(3k+3)/2}\}$,
$S_{b+2}(u_0)=\{u_{(3k+3)/2}\}$ and $S_i(u_0)=\emptyset$ for $i>b+2$;

\item[(iv)] if $k$ is even, then
$S_{b+1}(u_0)=\{u_{\pm((3k+2)/2)}\}$ and $S_i(u_0)=\emptyset$ for $i> b+1$.
\end{itemize}
\end{lemma}

\begin{proof}
By a careful inspection of the neighbors' sets of vertices
$u_i$ and $v_i$ (and using the assumption that $k\ge 13$),
we get that (i) holds
(see also Figure~\ref{GPG_struc}).

\begin{figure}[ht]
\begin{center}
 \includegraphics[width=0.70\hsize]{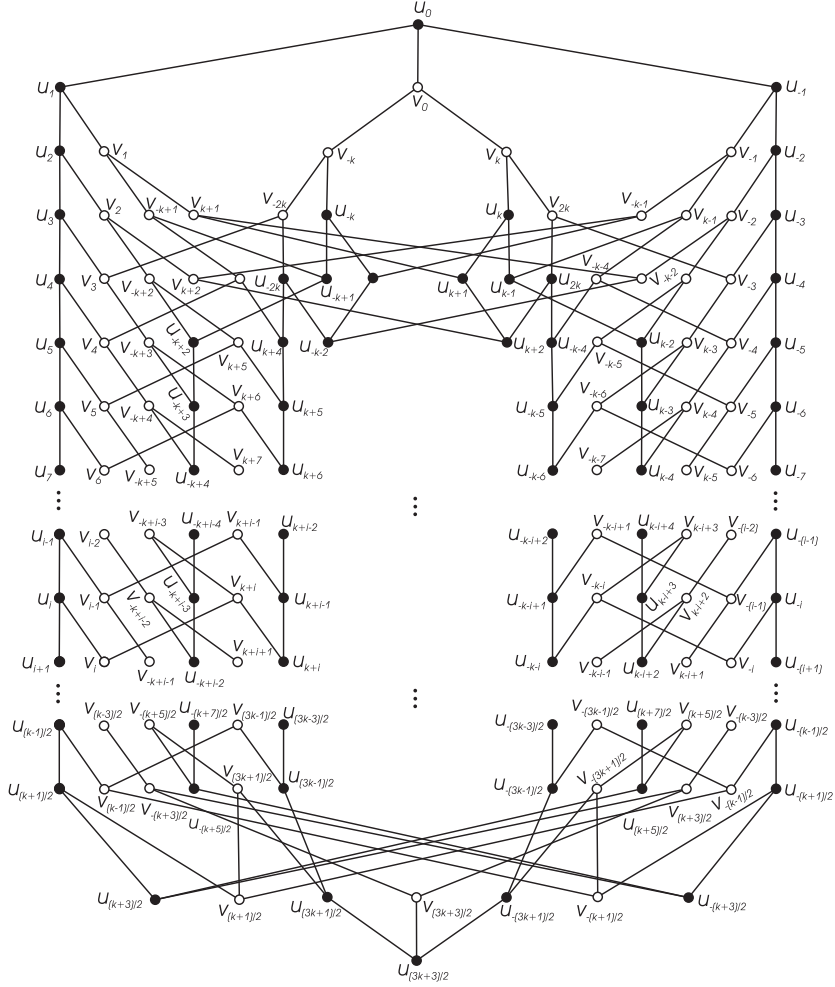}
\caption{\small The generalized Petersen graph $\GP(n,k)$ where
$k\ge 13$ is odd and $n=3k+3$.    \label{GPG_struc}}
\end{center}
\end{figure}

We now prove (ii) using induction. Similarly as in
the proof of (i) above we see that (ii) holds for $i \in \{6,7\}$.
Let us now assume that (ii) holds for $i-1$ and $i$,
where $i \in \{7,\dots,b-1\}$. Hence we have
$$
  S_{i-1}(u_0)=\{u_{\pm (i-1)}, u_{\pm(k+i-2)}, u_{\pm( k-i+4)},
  v_{\pm(i-2)}, v_{\pm( k-i+3)},v_{\pm(k+i-1)}\}
$$
and
$$
  S_i(u_0)=\{u_{\pm i}, u_{\pm(k+i-1)}, u_{\pm( k-i+3)},
  v_{\pm(i-1)}, v_{\pm( k-i+2)},v_{\pm(k+i)}\}.
$$
Obviously, $S_{i+1}(u_0)$ consists of all the neighbors of vertices in $S_i(u_0)$,
which are not in $S_{i-1}(u_0)$ or $S_{i}(u_0)$. Thus, by (\ref{genpet}),
$S_{i+1}(u_0)=\{u_{\pm (i+1)}, u_{\pm(k+i)}, u_{\pm(k-i+4)}, v_{\pm i}, v_{\pm (k-i+1)},v_{\pm(k+i+1)}\}$,
and the result follows (see also Figure~\ref{GPG_struc}).


\smallskip
Let us now prove (iii). If $k$ is odd, then $b=(k+1)/2$. By (ii),
$$
  S_{b-1}(u_0)=\{u_{\pm (k-1)/2}, u_{\pm(3k-3)/2}, u_{\pm(k+7)/2},
                 v_{\pm(k-3)/2}, v_{\pm(k+5)/2},v_{\pm(3k-1)/2}\},
$$
and
$$
  S_b(u_0)=\{u_{\pm(k+1)/2}, u_{\pm(3k-1)/2}, u_{\pm(k+5)/2},
             v_{\pm(k-1)/2}, v_{\pm(k+3)/2},v_{\pm(3k+1)/2} \}.
$$
Computing the neighbors of the vertices in $S_b(u_0)$ and sorting out those which are
in $S_{b-1}(u_0)$ or $S_b(u_0)$, we obtain
$S_{b+1}(u_0)=\{u_{\pm(k+3)/2},v_{\pm(k+1)/2}, u_{\pm(3k+1)/2}, v_{(3k+3)/2}\}$.
Furthermore, computing the neighbors of the vertices in $S_{b+1}(u_0)$ and sorting out
those which are in $S_b(u_0)$ or $S_{b+1}(u_0)$, we obtain  $S_{b+2}(u_0)=\{u_{(3k+3)/2}\}$.
Note that
$$
  \bigcup_{i=0}^{b+2}S_i(u_0)=V(\GP(n,k)),
$$
and hence the result follows.

\noindent
The proof of (iv) is similar to that of (iii) and is therefore left to the reader.
\end{proof}

We have the following immediate corollary of Lemma~\ref{lem:k}.

\begin{corollary}
\label{cor:k}
Let $k\ge 13$ be an integer, let $n=3k+3$, let $b=\lceil (k+1)/2\rceil$,
and let $u_0 \in V(\GP(n,k))$.
Then the following hold:
\begin{itemize}
\itemsep=0cm
\item[(i)] $|S_1(u_0)|=3 $, $|S_2(u_0)|=6  $, $|S_3(u_0)|=12$,
$|S_4(u_0)|=16$ and $|S_5(u_0)|=14$;

\item[(ii)] $|S_{i}(u_0)|=12$ for $6\le i\le b$;

\item[(iii)] if $k$ is odd, then
$|S_{b+1}(u_0)|=7$, $|S_{b+2}(u_0)|=1$ and $|S_i(u_0)|=0$ for $i>b+2$;

\item[(iv)] if $k$ is even, then $|S_{b+1}(u_0)|= 2$
and $|S_i(u_0)|= 0$ for $i> b+1$.
\end{itemize}
\end{corollary}

The proofs of the next lemma and corollary are omitted as they can
be carried out using the same arguments as in the proof of
Lemma~\ref{lem:k}. (Note that if $k \equiv -1 (\mod 3)$, then $2k+1$
is the multiplicative inverse of $k$ in $\ZZ_{3k+3}$.)

\begin{lemma}
\label{lem:k^-1}
Let $k\ge 13$ be an integer, let $n=3k+3$, let $b=\lceil (k+1)/2\rceil$,
and let $u_0 \in V(\GP(n,2k+1))$.
Then the following hold:
\begin{itemize}
\itemsep=0cm
\item[(i)] $S_1(u_0)=\{u_{\pm 1},v_0\}$, $S_2(u_0) =\{u_{\pm 2}, v_{\pm 1}, v_{\pm (k+2)}\}$, \\
          $S_3(u_0) =\{u_{\pm 3}, u_{\pm (k+2)}, v_{\pm 2}, v_{\pm (k+1)},v_{\pm (k+3)}, v_{\pm (k-1)}\}$, \\
          $S_4(u_0) =\{u_{\pm 4}, u_{\pm (k-1)}, u_{\pm (k+1)}, u_{\pm(2k)}, v_{\pm 3}, v_{\pm k},
                        v_{\pm (k+4)},  v_{\pm (k-2)}\}$, \\
          $S_5(u_0)=\{u_{\pm 5}, u_{\pm k}, u_{\pm (k+4)}, u_{\pm (k-2)}, v_{\pm 4}, v_{\pm (k+5)}, v_{\pm (k-3)}\} $;

\item[(ii)] $S_{i}(u_0)=\{u_{\pm i}, v_{\pm(i-1)},
u_{\pm(k+i-1)}, u_{\pm(2k+i)}, v_{\pm(2k+i+1)},v_{\pm(k+i)}\}$ for $6\le i\le b$;

\item[(iii)] if $k$ is odd, then $S_{b+1}(u_0)=\{u_{\pm(k+3)/2}, v_{\pm(k+1)/2}, u_{\pm(3k+1)/2}, v_{(3k+3)/2}\}$,\\
            $S_{b+2}(u_0)=\{u_{(3k+3)/2}\}$, and $S_i(u_0)=\emptyset$ for $i> b+2$;

\item[(iv)] if $k$ is even, then
$S_{b+1}(u_0)=\{u_{\pm(3k+2)/2}\}$, and $S_i(u_0)=\emptyset$ for $i> b+1$.
\end{itemize}
\end{lemma}

\begin{corollary}
\label{cor:k^-1}
Let $k\ge 13$ be an integer, let $n=3k+3$, let $b=\lceil (k+1)/2\rceil$,
and let $u_0 \in V(\GP(n,2k+1))$.
Then the following hold:

\begin{itemize}
\itemsep=0cm
\item[(i)] $|S_1(u_0)|=3 $, $|S_2(u_0)|=6$,
$|S_3(u_0)|=12$, $|S_4(u_0)|=16$, $|S_5(u_0)|=14$;

\item[(ii)] $|S_{i}(u_0)|=12$ for $6\le i\le b$;

\item[(iii)] if $k$ is odd, then
$|S_{b+1}(u_0)|=7$, $|S_{b+2}(u_0)|=1$ and $|S_i(u_0)|=0$ for $i>b+2$;

\item[(iv)] if $k$ is even, then $|S_{b+1}(u_0)|=
2$ and $|S_i(u_0)|=0$ for $i> b+1$.
\end{itemize}
\end{corollary}

\newpage

We are now ready to prove the main result of this section.

\begin{theorem}
\label{thm:glavni}
Let $k$ be a positive integer. Then the following hold:

\begin{itemize}
\itemsep=0cm
\item[(i)]  if $k\equiv 0 (\mod 3)$, then
$\GP(3k+3,k)$ is not strongly distance--balanced;

\item[(ii)]  if $k\not\equiv 0 (\mod 3)$, then
$\GP(3k+3,k)$ is strongly distance--balanced, or it is isomorphic
to $\GP(9,2)$, which is not strongly distance--balanced;
\item[(iii)] if $k \ge
2$ and $k\equiv 0 (\mod 3)$, then $\GP(3k-3,k)$ is not strongly
distance--balanced;
\item[(iv)] if $k \ge 2$ and
$k\not\equiv 0 (\mod 3)$, then $\GP(3k-3,k)$ is strongly
distance--balanced, or it is isomorphic to $\GP(9,4)\cong \GP(9,2)$, which is
not strongly distance--balanced.
\end{itemize}
\end{theorem}

\begin{proof}
Part (i) can be easily verified for $k \le 18$, so
assume that $k \ge 21$. Let us suppose that, by contradiction,
$\GP(3k+3,k)$ is strongly distance--balanced.

We distinguish two different cases depending on the parity of $k$.
Assume first that $k$ is odd.
By Lemma~\ref{lem:k}, the largest distance of some vertex from
$u_0$ is equal to $d=(k+5)/2$; in fact $S_{d}(u_0)=\{u_{(3k+3)/2}\}$.
Since $\GP(3k+3,k)$ is strongly distance--balanced,
it follows that $D^{d}_{d+1}(u_0,v_0)=\emptyset$.
Moreover, since $k$ is odd and $3k+3$ is even,
we have that $\GP(3k+3,k)$ is bipartite,
and hence $D^d_d(u_0,v_0)=\emptyset$.
Therefore $D^{d}_{d-1}(u_0,v_0)=\{u_{(3k+3)/2}\}$. Since
$\GP(3k+3,k)$ is strongly distance--balanced, it follows that
$|D^{d-1}_{d}(u_0,v_0)|=1$. Further, by (\ref{grupa}), we have $u_i\in
D^{d-1}_d(u_0,v_0)$ if and only if $ u_{-i}\in
D^{d-1}_d(u_0,v_0)$. Similarly, $v_i\in D^{d-1}_d(u_0,v_0)$ if and
only if $ v_{-i}\in D^{d-1}_d(u_0,v_0)$. It follows that
$D^{d-1}_{d}(u_0,v_0)=\{v_{(3k+3)/2}\}$. But the vertex
$v_{(3k+3)/2}$ belongs to the $(k+1)$-cycle
$$
(v_0,v_k,v_{2k},\ldots, v_{(3k+3)/2}, \ldots,v_0),
$$
and thus
$d(v_0,v_{(3k+3)/2})\le (k+1)/2=d-2$. This contradiction completes
the proof of (i) in the case when $k$ is odd.

\medskip

Assume next that $k$ is even and let $d=(k+4)/2$.
We first show that $D_i^i(u_0,v_0) = \emptyset$ for
$i \in \{1, \ldots, d-2\}$.
Suppose on contrary that $D_i^i(u_0,v_0) \ne \emptyset$ for some
$i \in \{1, \dots, d-2\}$, and let $j$ be the smallest positive integer such that
$D_j^j(u_0,v_0) \ne \emptyset$. Let $x \in D_j^j(u_0,v_0)$. Since $x$ is at distance
$j$ from $v_0$, we must have $S_1(x) \cap D_{j-1}^j(u_0,v_0) \ne \emptyset$
by minimality of $j$. Therefore there is an edge between two vertices
from $S_j(u_0)$. However, the sphere $S_j(u_0)$ is given in
Lemma~\ref{lem:k}, and it is easy to check that this is not possible.
Hence $D_i^i(u_0,v_0) = \emptyset$.

By Lemma~\ref{lem:k}, $S_d(u_0)=\{u_{3k/2+1}, u_{-3k/2-1}\}$.
Observe that $d(v_0,u_{3k/2+1}) \le d-1$, since $v_0, v_{2k+3},$
$v_{k+3}, u_{k+3}, u_{k+4}, \ldots, u_{3k/2+1}$ is a path of
length $d-1$ between $v_0$ and $u_{3k/2+1}$. Moreover, by the
triangle inequality, $d(v_0,u_{3k/2+1}) = d-1$. Similarly,
$d(v_0,u_{-3k/2-1}) = d-1$. Therefore $D_{d+1}^d(u_0,v_0)
= D_d^d(u_0,v_0) = \emptyset$ and
$D_{d-1}^d(u_0,v_0)=\{u_{3k/2+1}, u_{-3k/2-1}\}$. Combining together
Corollary~\ref{cor:k} and the fact that $\GP(3k+3,k)$ is strongly
distance--balanced, we can now compute the cardinalities of the sets
$D_{i-1}^i(u_0,v_0), D_i^{i-1}(u_0,v_0)$, $i \in \{1, \ldots, d)$
and $D_{d-1}^{d-1}(u_0,v_0)$. In particular,
$|D^{d-1}_{d-2}(u_0,v_0)|=6$, $|D_{d-1}^{d-1}(u_0,v_0)|=4$ and
$|D_d^{d-1}(u_0,v_0)|=2$.

\noindent
Observe that, by Lemma~\ref{lem:k}, we have
$$
  D_{d-2}^{d-1}(u_0,v_0) \cup D_{d-1}^{d-1}(u_0,v_0) \cup D_d^{d-1}(u_0,v_0) =
$$
$$
  \{ u_{\pm (k/2+1)}, v_{\pm k/2}, v_{\pm (k/2+1)}, u_{\pm (k/2+2)}, u_{\pm 3k/2},
     v_{\pm(3k/2+1)} \}.
$$
Since the vertices $v_{k/2}$ and $v_{-k/2}$ are contained on the cycle
\begin{equation}
\label{cikel}
C = (v_0, v_k, v_{2k}, \ldots, v_{0})
\end{equation}
of length $k+1$, we have
$d(v_0,v_{k/2}) \le k/2=d-2$ and $d(v_0,v_{-k/2}) \le k/2=d-2$. Hence,
$v_{k/2}, v_{-k/2} \in D_{d-2}^{d-1}(u_0,v_0)$. Furthermore, the path
$v_0, v_k, u_k, u_{k-1}, \ldots, u_{k/2+2}$ has length $k/2=d-2$, implying
$u_{k/2+2} \in D_{d-2}^{d-1}(u_0,v_0)$. Similarly we get that
$u_{-(k/2+2)} \in D_{d-2}^{d-1}(u_0,v_0)$.
Finally, since the vertices $v_{3k/2}$ and $v_{-3k/2}$ are also contained on the cycle C
in (\ref{cikel}) above, and since $d(v_0,v_{k/2})=d(v_0,v_{-k/2})=d-2$, we must have
that $d(v_0,v_{3k/2}) \le d-3$ and $d(v_0,v_{-3k/2}) \le d-3$. But this now implies
$d(v_0,u_{3k/2}) \le d-2$ and $d(v_0,u_{-3k/2}) \le d-2$, and hence
$u_{3k/2}, u_{-3k/2} \in D_{d-2}^{d-1}(u_0,v_0)$.
Since $u_{k/2+1}$ is adjacent with $u_{k/2+2}$ we have
$u_{k/2+1} \in D_{d-1}^{d-1}(u_0,v_0)$. Similarly, we get
$u_{-k/2-1} \in D_{d-1}^{d-1}(u_0,v_0)$.

We now show $v_{k/2+1} \in D_d^{d-1}(u_0,v_0)$. Suppose $v_{k/2+1}
\in D_{d-1}^{d-1}(u_0,v_0)$. Then, by (\ref{grupa}), $v_{-k/2-1}
\in$ $D_{d-1}^{d-1}(u_0,v_0)$. Since $\GP(3k+3,k)$ is strongly
distance--balanced, we have that $v_{3k/2+1}, v_{-3k/2-1} \in
D_d^{d-1}(u_0,v_0)$. Furthermore, since
$D_{d-2}^{d-2}(u_0,v_0)=\emptyset$, we get $S_1(v_{k/2+1}) \cap
D_{d-1}^{d-2}(u_0,v_0) \ne \emptyset$ and $S_1(v_{k/2+1}) \cap
D_{d-2}^{d-1}(u_0,v_0) \ne \emptyset$. But this is impossible since
$u_{k/2+1} \in D_{d-1}^{d-1}(u_0,v_0)$ and $v_{3k/2+1} \in
D_d^{d-1}(u_0,v_0)$.

In a similar fashion we can show that $v_{3k/2+1} \in D_d^{d-1}(u_0,v_0)$.
But then, by (\ref{grupa}), we have that also
$v_{-k/2-1}, v_{-3k/2-1} \in D_d^{d-1}(u_0,v_0)$, a contradiction.
This completes the proof of part (i).

\medskip \noindent
To prove part (ii) suppose first that $k\equiv 1 (\mod 3)$. Then it
is easy to check that $k^2\equiv 1 (\mod 3k+3)$. Hence $\GP(3k+3,k)$
is vertex--transitive and, by Corollary~\ref{cor:vertex}, strongly
distance--balanced.

Suppose next that $k\equiv -1 (\mod 3)$. For $k\le 12$, we have
verified the strongly distance--balanced property of generalized
Petersen graphs $\GP(3k+3,k)$ with program package {\sc Magma}
\cite{Magma}. In particular, $\GP(9,2)$ is the only graph among the
generalized Petersen graphs $\GP(9,2)$, $\GP(18,5)$, $\GP(27,8)$ and
$\GP(36,11)$, which is not strongly distance--balanced. We may
therefore assume that $k\ge 13$. Observe that $3k+3$ and $k$ are
relatively prime and that $k(2k+1) \equiv 1 (\mod 3k+3)$. Hence, by
(\ref{iso}), $\GP(3k+3,k) \cong \GP(3k+3,2k+1)$. Combining together
Corollaries~\ref{cor:k} and \ref{cor:k^-1}, we get
$|S_i(u_0)|=|S_i(v_0)|$ for all integers $i$. Finally, by
(\ref{grupa}), we have also that $|S_i(u_0)|=|S_i(u_t)|$ and
$|S_i(v_0)|=|S_i(v_t)|$ for all integers $i$ and for all $t \in
\ZZ_n$, completing the proof of part (ii).

\medskip \noindent
The proof of part (iii) is analogous to the proof of part
(i) and is thus omitted.

\medskip \noindent
Finally, to prove part (iv),
assume first $k\equiv  -1 (\mod 3)$. Then it is easy to check that
$k^2\equiv 1(\mod 3k-3)$. Thus $\GP(3k-3,k)$ is vertex--transitive, and so, by
Corollary~\ref{cor:vertex}, strongly distance--balan\-ced.

Next assume $k \equiv 1 (\mod 3)$. Observe that in this case the
multiplicative inverse of $k$ in $\ZZ_{3k-3}$ is $2k-1$. Hence
$\GP(3k-3,k) \cong \GP(3k-3,2k-1)$ by~(\ref{iso}). Furthermore, we
have $\GP(3k-3,2k-1) \cong \GP(3k-3,(3k-3)-(2k-1)) =
\GP(3(k-2)+3,k-2)$. But then part (ii) implies that the graph $\GP(3k-3,k)$ is
either strongly distance--balanced or isomorphic to $\GP(9,4) \cong GP(9,2)$,
as required.
\end{proof}

To wrap up this section, let us remark that
an application of similar methods to the ones used in the proof of
Theorem~\ref{thm:glavni}, leads us to the following result
identifying another infinite family of strongly distance--balanced
generalized Petersen graphs.

\begin{theorem}
\label{thm:2k+2}
Let $k$ be a positive integer. Then $\GP(2k+2,k)$
is strongly distance--balanced if and only if $k$ is odd.
\end{theorem}

%


\end{document}